\documentclass{amsart}
\usepackage{amsmath}
\usepackage{amsfonts}
\usepackage{array}
\usepackage{amsthm}
\makeatletter
\def\thm@space@setup{%
  \thm@preskip=\parskip \thm@postskip=0pt
}
\makeatother
\usepackage{amssymb}
\usepackage{parskip}
\newcommand{\sr}{\mathbb{R}}

\newtheorem{lemma}
{Lemma}
\newtheorem{cor}
{Corollary}
\newtheorem{prop}
{Proposition}
{Conjecture}
{Question}
{Sublemma}
\newtheorem{claim}
{Claim}
\newtheorem{theorem}{Theorem}

\theoremstyle{definition}
{Definition}

\theoremstyle{remark}

\newtheorem*{acknowledgement}{Acknowledgement}

\numberwithin{equation}{section}

\begin{document}

\title[A sharp bilinear estimate for the Klein--Gordon equation]{A sharp bilinear estimate for the Klein--Gordon equation in
arbitrary space-time dimensions}
\author{Chris Jeavons}

\address{School of Mathematics, The Watson Building, The University of Birmingham, Edgbaston, Birmingham, B15 2TT, England}
\email{jeavonsc@maths.bham.ac.uk}

\maketitle
\begin{abstract}
We prove a sharp bilinear inequality for the Klein--Gordon equation
on $\sr^{d+1}$, for any $d \geq 2$. This extends work of
Ozawa--Rogers and Quilodr\'an for the Klein--Gordon equation and
generalises work of Bez--Rogers for the wave equation. As a
consequence we obtain a sharp Strichartz estimate for the solution
of the Klein--Gordon equation in five spatial dimensions for data
belonging to $H^1$. We show that maximisers for this estimate do not
exist and that any maximising sequence of initial data concentrates
at spatial infinity.
\end{abstract}

\section{Introduction}
For the Klein--Gordon equation on $\sr^{1+1}$, very recently in
\cite{O--R'?} it was shown that the bilinear estimate
\begin{align}\label{R--O '12}
\left\|e^{it \sqrt{1-\Delta}}f_1\right. & \left. e^{it \sqrt{1-\Delta}}f_2\right\|_{L^2(\sr^{2})}^2 \nonumber \\
& \leq\frac{1}{(2\pi)^2}\int_{\sr^{2}}\left|\smash{\widehat{f_1}(y_1)}\right|^2 \left|\smash{\widehat{f_2}(y_2)}\right|^2\frac{(1+y_1^2)^{\frac{3}{4}}(1+y_2^2)^{\frac{3}{4}}}{\left|y_1-y_2\right|}\, \mathrm{d}y_1 \mathrm{d}y_2
\end{align}
holds whenever $f_1$ and $f_2$ have disjoint Fourier supports, and
that the constant $\frac{1}{(2\pi)^2}$ is sharp. The main motivation
behind the present paper was to identify a natural generalisation of
this sharp bilinear estimate to arbitrary dimensions. In achieving
this, we simultaneously extend work of Quilodr\'an in \cite{Q'12}
and generalise work of Bez--Rogers \cite{B--R'11}. We will also
obtain a new Strichartz estimate with sharp constant for the
Klein--Gordon equation on $\mathbb{R}^{5+1}$ with $H^1$-initial
data.

Throughout this paper, we let $\, \widehat{.} \,$ denote the spatial
Fourier transform on $\sr^d$, defined on the Schwartz class as
\begin{equation*}
\widehat{f}(\xi)=\int_{\sr^d}f(x)e^{-ix\cdot\xi}\,\mathrm{d}x, \;\;\; \xi\in\sr^d.
\end{equation*}
For fixed $s \geq 0$, define also the Klein--Gordon propagator
$e^{it\phi_s(\sqrt{-\Delta})}$ by
\begin{equation} \label{KG-def}
e^{it\phi_s(\sqrt{-\Delta})}f(x)=\frac{1}{(2\pi)^d}\int_{\sr^d}\widehat{f}(\xi)e^{ix\cdot\xi+it(s^2+\left|\xi\right|^2)^{\frac{1}{2}}}\,\mathrm{d}\xi, \;\; x\in\sr^d,\;\; t\in\sr,
\end{equation}
on the Schwartz class, where $\phi_s(r)=\sqrt{s^2+r^2}$, for
$r\in\sr$. In order to state the main result, in what follows we let
\begin{equation*}
K_s(y_1, y_2)=\frac{\left(\phi_s(\left|y_1\right|)\phi_s(\left|y_2\right|)-y_1\cdot y_2-s^2\right)^{\frac{d-2}{2}}}{\left(\phi_s(\left|y_1\right|)\phi_s(\left|y_2\right|)-y_1\cdot y_2+s^2\right)^{\frac{1}{2}}}
\end{equation*}
be a function on $\sr^{2d}$, and we introduce the constant
\begin{equation*}
\mathbf{KG}(d)=\frac{2^{-\frac{d-1}{2}}\left|\mathbb{S}^{d-1}\right|}{(2\pi)^{3d-1}}
\end{equation*}
for $d\geq 1$, which will appear throughout the paper.
\begin{theorem}\label{bilineart}
If $d\geq 2$ and $s\geq 0$, then
\begin{align}\label{bilinear-eq}
& \left\|e^{it\phi_s(\sqrt{-\Delta})} f_1 e^{it\phi_s(\sqrt{-\Delta})}f_2 \right\|_{L^2(\sr^{d+1})}^2 \nonumber \\
& \qquad \qquad \leq \mathbf{KG}(d)\int_{\sr^{2d}}\left|\smash{\widehat{f_1}}(y_1)\right|^2\left|\smash{\widehat{f_2}}(y_2)\right|^2 \phi_s(\left|y_1\right|)\phi_s(\left|y_2\right|)K_s(y_1,y_2)\,\mathrm{d}y_1 \,\mathrm{d}y_2,
\end{align}
where the constant $\mathbf{KG}(d)$ is best possible since we have
equality for functions of the form
\begin{equation} \label{maxdef}
\widehat{f_1}(\xi)=\widehat{f_2}(\xi)=\frac{e^{-a\phi_s(\left|\xi\right|)
}}{\phi_s(\left|\xi\right|)},
\end{equation}
for $a 
>0$. 
Further, if $d=1$, $s>0$ and $\widehat{f_1}, \widehat{f_2}$ have
disjoint support, or $s=0$ and $\widehat{f_1}, \widehat{f_2}$ have
disjoint angular support, then we have
\begin{align}\label{main-identity}
&\left\|e^{it\phi_s(\sqrt{-\Delta})}f_1 e^{it\phi_s(\sqrt{-\Delta})}f_2 \right\|_{L^2(\sr^{1+1})}^2 \nonumber \\
& \qquad \qquad = \frac{\mathbf{KG}(1)}{2} \int_{\sr^{2d}}\left|\smash{\widehat{f_1}(y_1)}\right|^2\left|\smash{\widehat{f_2}(y_2)}\right|^2 \phi_s(\left|y_1\right|)\phi_s(\left|y_2\right|)K_s(y_1,y_2)\,\mathrm{d}y_1 \,\mathrm{d}y_2.
\end{align}
\end{theorem}
It turns out that the functions described in \eqref{maxdef} above
play an important role in some of the applications of this
inequality, as we will see below in Corollary \ref{cor-1c}.

In the case $d=1$, we observe that
\begin{equation}\label{pointwise}
K_s(y_1 ,y_2)\leq \frac{1}{s^2}\frac{(s^2+y_1^2)^{\frac{1}{4}}(s^2+y_2)^{\frac{1}{4}}}{\left|y_1-y_2\right|}
\end{equation}
for almost every $(y_1,y_2)\in\sr^2$ with $y_1\neq y_2$ and $s>0$.
One can see this by first reducing to the case $s=1$ and then a
direct argument shows that the claimed inequality is equivalent to
\begin{equation*}
0\leq (y_1^2+y_2^2)(y_1-y_2)^4+y_1^2 y_2^2(y_1-y_2)^4,
\end{equation*}
which is clearly true. Since $\frac{1}{2}\mathbf{KG}(1) =
\frac{1}{(2\pi)^2}$, we see that \eqref{R--O '12} follows from
\eqref{main-identity}\footnote{in fact, the argument in
\cite{O--R'?} leading to \eqref{R--O '12} goes via the identity
\eqref{main-identity}, and they prove \eqref{pointwise} differently
using some trigonometric identities} and we claim that
\eqref{bilinear-eq} in Theorem \ref{bilineart} provides a natural
generalisation of this to higher dimensions.

Furthermore, one can deduce certain Strichartz estimates from
\eqref{bilinear-eq} with sharp constants, some of which recover
sharp Strichartz estimates due to Quildodr\'an in \cite{Q'12} and
Bez--Rogers in \cite{B--R'11}, and we also obtain a new sharp
Strichartz estimate for the Klein--Gordon equation in five spatial
dimensions (see the forthcoming Corollary \ref{cor-1c}). In order to
describe these results, we define the inhomogeneous Sobolev norm as
follows
\begin{equation*}
\left\|f\right\|_{H^m(\sr^d)}=\left\|\smash{(1+\left|\,.\,\right|^2)^{\frac{m}{2}}\widehat{f}}\right\|_{L^2(\sr^d)}, \, m\geq 0.
\end{equation*}

When $d=2$ and $s=1$, we have
\begin{equation}\label{one}
K_s(y_1, y_2)= (\phi_1(\left|y_1\right|)\phi_1(\left|y_2\right|)-y_1\cdot y_2+1)^{-\frac{1}{2}} \leq 2^{-\frac{1}{2}},
\end{equation}
for $(y_1, y_2)\in\sr^4$. Taking $f_1=f_2$ in \eqref{bilinear-eq} it
follows that
\begin{equation} \label{e:Q24}
\| e^{it\sqrt{1-\Delta}} f\|_{L^4(\sr^{2+1})} \leq \frac{1}{2^{5/4}\pi} \|f\|_{H^{1/2}(\mathbb{R}^2)}.
\end{equation}
Estimate \eqref{e:Q24} is due to Quilodr\'an \cite{Q'12} and he
showed that the constant is sharp but that maximisers do not
exist\footnote{in \cite{Q'12} the perspective is that of adjoint
Fourier restriction inequalities for the hyperboloid, and we choose
to present the estimates from \eqref{e:Q24} in terms of the
Klein--Gordon propagator}.

Similarly, when $d=3$ and $s=1$, we get
\begin{equation}\label{two}
K_s(y_1, y_2)=\frac{(\phi_1(\left|y_1\right|)\phi_1(\left|y_2\right|)-y_1\cdot y_2-1)^{\frac{1}{2}}}{(\phi_1(\left|y_1\right|)\phi_1(\left|y_2\right|)-y_1\cdot y_2+1)^{\frac{1}{2}}}\leq 1,
\end{equation}
and \eqref{bilinear-eq} implies
\begin{equation} \label{e:Q34}
\| e^{it\sqrt{1-\Delta}} f\|_{L^4(\sr^{3+1})} \leq \frac{1}{(2\pi)^{7/4}} \|f\|_{H^{1/2}(\mathbb{R}^3)}.
\end{equation}
Again, the constant is sharp and maximisers do not exist (due to
Quilodr\'an \cite{Q'12}). We remark that we prove Theorem
\ref{bilineart} using the approach of Foschi in \cite{F'06}, as did
Quilodr\'an, and so it is not at all a surprise that \eqref{e:Q24}
and \eqref{e:Q34} follow from Theorem \ref{bilineart}.

In this paper, we obtain the following new sharp form of a classical
Strichartz estimate for the full solution of the Klein--Gordon
equation for data in the energy space.
\begin{cor}\label{cor-2c}
Suppose that 
$\partial_{tt} u-\Delta u+u=0$ on $\sr^{5+1}$, 
then
\begin{align}\label{cor-1.2}
\left\|u\right\|_{L^4(\sr^{5+1})}
\leq & \left(\frac{1}{8\pi}\right)^{\frac{1}{2}}\left(\left\|u(0)\right\|_{H^1(\sr^5)}^2+\left\|\partial_tu(0)\right\|_{L^2(\sr^5)}^2\right)^{\frac{1}{2}} 
.
\end{align}
The constant $\left(\frac{1}{8\pi}\right)^{\frac{1}{2}}$ is sharp,
but there are no nontrivial functions for which we have equality.
\end{cor}
A nonsharp form of \eqref{cor-1.2} was proved by Strichartz in
\cite{S'77}. The sharp inequality \eqref{cor-1.2} is deduced from
the following sharp estimate for the one-sided propagator
$e^{it\phi_s(\sqrt{-\Delta})}$. In order to state this result, we introduce the notation
\[\left\|f\right\|_{(s)}^4:=\left\|\phi_s(\sqrt{-\Delta})f\right\|_{L^2}^4-s^2\left\|(\phi_s(\sqrt{-\Delta}))^{\frac{1}{2}}f\right\|_{L^2}^4.\]
Notice that if $s=1$ then $\left\|f\right\|_{(s)}$ may be bounded
above by the inhomogeneous norm $\left\|.\right\|_{H^1}$, and if
$s=0$ then $\left\|f\right\|_{(s)}$ is just the $\dot{H}^1$-norm of
$f$.
\begin{cor}\label{cor-1c}
Let $s \geq 0$. Then
\begin{equation}\label{cor-1.1}
\left\|e^{it\phi_s(\sqrt{-\Delta})}f\right\|_{L^4(\sr^{5+1})}
\leq \left(\frac{1}{24\pi^2}\right)^{\frac{1}{4}}\left\|f\right\|_{(s)}.
\end{equation}
The constant $\left(\frac{1}{24\pi^2}\right)^{\frac{1}{4}}$ is sharp as we have the maximising sequence
$\left(g_a\right)_{a>0}$ defined by
\begin{equation}\label{g_a}
g_a=\frac{f_a}{\left\|f_a\right\|_{(s)}},
\end{equation}
where
\begin{equation} \label{f_a}
\widehat{f_a}(\xi)=\frac{e^{-a\phi_s(\left|\xi\right|)}}{\phi_s(\left|\xi\right|)}
\end{equation}
as $a\rightarrow 0+$
, but when $s>0$ there are no functions for which we have equality.
\end{cor}
By maximising sequence for \eqref{cor-1.1} we mean a sequence of functions $(g_n)_{n\geq 1}$ satisfying $\left\|g_n\right\|_{(s)}\leq 1$ for which
\begin{equation*}
\left\|e^{it\phi_s(\sqrt{-\Delta})}g_n\right\|_{L^4(\sr^{5+1})}\rightarrow \left(\frac{1}{24\pi^2}\right)^{\frac{1}{4}}
\end{equation*}
as $n\rightarrow\infty$.

When $s=1$, \eqref{cor-1.1} is the sharp estimate
\begin{equation} \label{e:J54}
\| e^{it\sqrt{1-\Delta}} f\|_{L^4(\sr^{5+1})} \leq \left(\frac{1}{24\pi^2}\right)^{\frac{1}{4}}(\|f\|_{H^{1}(\mathbb{R}^5)}^4 - \|f\|_{H^{1/2}(\mathbb{R}^5)}^4)^{1/4},
\end{equation}
which is a refinement of the sharp estimate
\begin{equation}\label{l4h1}
\left\|\smash{e^{it\sqrt{1-\Delta}}f}\right\|_{L^4(\sr^{5+1})}\leq \left(\frac{1}{24\pi^2}\right)^{\frac{1}{4}}\left\|f\right\|_{H^1(\sr^5)}.
\end{equation}
Both \eqref{e:J54} and \eqref{l4h1} are new. With nonsharp constant,
\eqref{l4h1} follows from \cite{S'77}.
That the constant in \eqref{l4h1} is sharp follows from the
observation that for the functions $f_a$ defined by \eqref{f_a} one has that
\begin{equation}\label{h1/2-decay}
a^5\left\|\phi_s(\sqrt{-\Delta})^{\frac{1}{2}}f_a\right\|^2_{L^2(\sr^5)}\rightarrow 0
\end{equation}
as $a\rightarrow 0+$ (see Section \ref{section:corollaries}). In
fact, a similar property holds for maximising sequences for \eqref{cor-1.1},
as we will see in our forthcoming Proposition \ref{conc}.

At this point, we make some remarks concerning the particular case
of the wave equation corresponding to the case $s=0$. With the
emphasis not on sharp constants, Klainerman--Machedon first
established bilinear estimates in the spirit of \eqref{bilinear-eq}
with different kinds of weights in the case $s=0$ (see \cite{KM'93},
\cite{KM'96}, \cite{KM'97}). Regarding sharp estimates, Theorem
\ref{bilineart} and \eqref{cor-1.1} for $s=0$ were established in
\cite{B--R'11} (see also \cite{C'09} for similar results for the
Schr\"odinger propagator). Also, maximisers exist in both cases
$s=0$ and $s > 0$ in Theorem \ref{bilineart}. However, in Corollary
\ref{cor-1c}, it is true that when $s=0$ and for any $a > 0$, the
function $f_a$ given by \eqref{f_a} is a maximiser, but when $s >
0$, there are no maximisers and when suitably normalised, the functions $f_a$ form a
maximising sequence as $a$ tends to zero.

As our final main result in this paper, we establish that \emph{any}
maximising sequence for the estimate \eqref{cor-1.1} must concentrate
at spatial infinity in the following precise sense.
\begin{prop}\label{conc} If $\left(g_n\right)_{n\geq 1}$ is any maximising sequence for \eqref{cor-1.1}, then for each $\varepsilon, R>0$ there exists $N\in\mathbb{N}$ so that if $n\geq N$,
\[\left\|\phi_s(\sqrt{-\Delta})^{\frac{1}{2}}g_n\right\|_{L^2(\sr^5)}<\varepsilon,\]
and
\begin{equation}\label{conc-ineq}
\left\|\widehat{\phi_s(\sqrt{-\Delta})g_n}\right\|_{L^2(B(0,R))} < \varepsilon,
\end{equation}
where $B(0,R)$ denotes the ball of radius $R$ centered at the origin
in $\sr^5$.
\end{prop}
The motivation for this result comes from the observation that the
particular maximising sequence $(g_a)$ considered in Corollary
\ref{cor-1c} satisfies these conditions. A result analogous to
\eqref{conc-ineq} was established in \cite{Q'12}, where it was shown
that any maximising sequence for either \eqref{e:Q24} or
\eqref{e:Q34} must concentrate at spatial infinity. We also remark
here that in the case $s=1$, Proposition \ref{conc} may be
interpreted as a statement about the concentration of the $H^1$-norm
of a maximising sequence for the inequality \eqref{cor-1.1}.

Largely as a result of the influential work of Foschi \cite{F'06}, a
body of very recent work has emerged on sharp constants and the
existence or nature of maximisers for space-time estimates
associated with dispersive PDE, to which this work belongs. In
addition to \cite{B--R'11}, \cite{F'06}, \cite{O--R'?} and
\cite{Q'12} already mentioned, see for example, \cite{BBCH'09},
\cite{BS'?}, \cite{C'09}, \cite{DMR'11}, \cite{HZ'06}, \cite{S'92},
and \cite{W'02} for sharp constants, and \cite{B'10}, \cite{CS1'12},
\cite{CS2'12}, \cite{FVV'11}, \cite{FVV'12}, \cite{K'03},
\cite{R'12}, and \cite{S'09} for results on maximisers.

\section{Proof of Theorem 1}
\subsection{The case $d\geq 2$ and $s\geq 0$}
In this section we will use the space-time Fourier transform,
defined for suitable functions $f$ on $\sr^{d+1}$ by
\begin{equation*}
\widetilde{f}(\xi,\tau)=\int_{\sr^{d+1}}f(x,t)e^{-i(t\tau+x\cdot\xi)}\,\mathrm{d}x\,\mathrm{d}t.
\end{equation*}
We note firstly that the space-time Fourier transform of
$v_j=e^{it\phi_s(\sqrt{-\Delta})}f_j$ will be the measure
\begin{equation*}
\widetilde{v_j}(\xi,\tau)=2\pi\delta(\tau-\phi_s(\left|\xi\right|))\widehat{f_j}(\xi)
\end{equation*}
for $j=1,2$, each supported on the hyperboloid in $\sr^{d+1}$,
\begin{equation*}
\left\{(y,(s^2+\left|y\right|^2)^{\frac{1}{2}}):y\in\sr^d\right\}.
\end{equation*}
Note that if $s>0$ and $d\geq 2$, the function defined by $K_s$ is
well-defined for any $y=(y_1, y_2)\in\sr^{2d}$. For example, if
$d=2$ the kernel reduces to
\begin{equation*}
K_s(y_1,y_2)=\frac{1}{(s^2+\phi_s(\left|y_1\right|)\phi_s(\left|y_2\right|)-y_1 \cdot y_2)^{\frac{1}{2}}},
\end{equation*}
and the denominator is always positive since
\begin{equation*}
y_1\cdot y_2\leq\left|y_1\right|\left|y_2\right|<(\left|y_1\right|^2+s^2)^{\frac{1}{2}}(\left|y_2\right|^2+s^2)^{\frac{1}{2}}+s^2=\phi_s(\left|y_1\right|)\phi_s(\left|y_2\right|)+s^2;
\end{equation*}
the claim for $d>2$ follows from this as the power $\frac{d-2}{2}$
is positive in this case.

If we now write $u_2=e^{it\phi_s(\sqrt{-\Delta})}f_1
e^{it\phi_s(\sqrt{-\Delta})}f_2$, then the space-time Fourier
transform of $u_2$ will be the convolution of the measures
$\widetilde{v_1}$ and $\widetilde{v_2}$, which may be written as
\begin{equation}\label{conv-form}
\displaystyle\widetilde{u_2}(\xi,\tau)=\frac{1}{(2\pi)^{d-1}}\int_{\sr^{2d}}\frac{\widehat{F}(y)}{(s^2+\left|y_1\right|^2)^{\frac{1}{4}}(s^2+\left|y_2\right|^2)^{\frac{1}{4}}}\delta\binom{\tau-\phi_s(\left|y_1\right|)-\phi_s(\left|y_2\right|)}{\xi-y_1-y_2}\,\mathrm{d}y,
\end{equation}
where $\xi\in\sr^d$ and $\tau\in\sr$ are fixed, we set
\begin{equation*}
\widehat{F}(y)=\widehat{f_1}(y_1)\widehat{f_2}(y_2)(s^2+\left|y_1\right|^2)^{\frac{1}{4}}(s^2+\left|y_2\right|^2)^{\frac{1}{4}},
\end{equation*}
and we use the notation $\delta\binom{t}{x}$ for the product
$\delta(t)\delta(x)$ on $\sr^{d+1}$. It is proved in \cite{Q'12}
that the function $\widetilde{u_2}$ is supported on the set
\begin{equation*}
\mathcal{H}_s=\left\{(\xi,\tau)\in\sr^{d+1}:\tau\geq((2s)^2+\left|\xi\right|^2)^{\frac{1}{2}}\right\},
\end{equation*}
for completeness we include the proof here. If $\xi=y_1+y_2$ and
$\tau=\phi_s(\left|y_1\right|)+\phi_s(\left|y_2\right|)$ we have
that
\begin{align*}
\tau^2 & =2s^2+\left|y_1\right|^2+\left|y_2\right|^2+2\phi_s(\left|y_1\right|)\phi_s(\left|y_2\right|) \\
& \geq \left|y_1\right|^2+\left|y_2\right|^2+2\left|y_1\right| \left|y_2\right|+4s^2 \\
& \geq 4s^2+\left|\xi\right|^2,
\end{align*}
since
\begin{equation}\label{ineq-1}
\phi_s(\left|y_1\right|)\phi_s(\left|y_2\right|)\geq s^2+\left|y_1\right|\left|y_2\right|,
\end{equation}
as can easily be seen by squaring both sides. By the Cauchy--Schwarz
inequality,
\begin{equation} \label{KG-CS}
\left|\widetilde{u_2}(\xi,\tau)\right|^2\leq\frac{I_s(\xi,\tau)}{(2\pi)^{2d-2}}\int_{\sr^{2d}}\left|\smash{\widehat{F}(y)}\right|^2K_s(y)\delta\binom{\tau-\phi_s(\left|y_1\right|)-\phi_s(\left|y_2\right|)}{\xi-y_1-y_2}\,\mathrm{d}y,
\end{equation}
where
\begin{equation*}
I_s(\xi,\tau)=\int_{\sr^{2d}}\frac{1}{K_s(y)\phi_s(\left|y_1\right|)\phi_s(\left|y_2\right|)}\delta\binom{\tau-\phi_s(\left|y_1\right|)-\phi_s(\left|y_2\right|)}{\xi-y_1-y_2}\,\mathrm{d}y.
\end{equation*}
Now, on the support of the delta measures, by the choice of $K_s$ we
have that
\begin{equation*}
K_s(y)=\frac{1}{2^{\frac{d-3}{2}}}\frac{(\tau^2-\left|\xi\right|^2-4s^2)^{\frac{d-2}{2}}}{(\tau^2-\left|\xi\right|^2)^{\frac{1}{2}}},
\end{equation*}
so that
\begin{align*}
I_s(\xi, \tau) & =\frac{2^{\frac{d-3}{2}}(\tau^2-\left|\xi\right|^2)^{\frac{1}{2}}}{(\tau^2-\left|\xi\right|^2-4s^2)^{\frac{d-2}{2}}}\int_{\sr^{2d}}\frac{1}{\phi_s(\left|y_1\right|)\phi_s(\left|y_2\right|)}\delta\binom{\tau-\phi_s(\left|y_1\right|)-\phi_s(\left|y_2\right|)}{\xi-y_1-y_2}\,\mathrm{d}y \\
& =\frac{2^{\frac{d-3}{2}}(\tau^2-\left|\xi\right|^2)^{\frac{1}{2}}}{(\tau^2-\left|\xi\right|^2-4s^2)^{\frac{d-2}{2}}}\sigma_s\ast\sigma_s(\xi,\tau),
\end{align*}
where we have defined the measure $\sigma_s$ on $\sr^{d+1}$ by
\begin{align*}
\int_{\sr^{d+1}}g(x,t)\,\mathrm{d}\sigma_s(x,t) & =\int_{\sr^{d+1}}g(x,t)\delta(t-\phi_s(\left|x\right|))\frac{\,\mathrm{d}x\,\mathrm{d}t}{\phi_s(\left|x\right|)}.
\end{align*}
Indeed, since $\xi=y_1+y_2$ we have that
$\left|\xi\right|^2=\left|y_1\right|^2+\left|y_2\right|^2+2y_1\cdot
y_2$, and since
$\tau=\phi_s(\left|y_1\right|)+\phi_s(\left|y_2\right|)$ it follows
that
\begin{equation*}
\tau^2=2s^2+\left|y_1\right|^2+\left|y_2\right|^2+2(s^2+\left|y_1\right|^2)^{\frac{1}{2}}(s^2+\left|y_2\right|^2)^{\frac{1}{2}},
\end{equation*}
and so we obtain
\begin{equation*}
\tau^2-\left|\xi\right|^2=2(s^2+\left|y_1\right|^2)^{\frac{1}{2}}(s^2+\left|y_2\right|^2)^{\frac{1}{2}}+2s^2-2y_1\cdot y_2,
\end{equation*}
so that
\begin{equation*}
\frac{1}{2}(\tau^2-\left|\xi\right|^2)=\phi_s(\left|y_1\right|)\phi_s(\left|y_2\right|)+s^2-y_1\cdot y_2,
\end{equation*}
and
\begin{equation*}
\phi_s(\left|y_1\right|)\phi_s(\left|y_2\right|)-s^2-y_1\cdot y_2=\frac{1}{2}(\tau^2-\left|\xi\right|^2-4s^2).
\end{equation*}
Hence we need to compute the quantity
\begin{equation*}
J_s(\xi,\tau):=\int_{\sr^{2d}}\frac{1}{\phi_s(\left|y_1\right|)\phi_s(\left|y_2\right|)}\delta\binom{\tau-\phi_s(\left|y_1\right|)-\phi_s(\left|y_2\right|)}{\xi-y_1-y_2}\,\mathrm{d}y=\sigma_s\ast\sigma_s(\xi,\tau).
\end{equation*}
It is known, \cite{S'77}, that the measure $\sigma_s$ is invariant
under Lorentz transformations, and hence so is the convolution
$J_s$. Using this invariance, the convolution may be computed
easily.
\begin{lemma}\label{measure-conv} For all $(\xi,\tau)\in\mathcal{H}_s$ we have that
\begin{equation*}
J_s(\xi,\tau)=\frac{\left|\mathbb{S}^{d-1}\right|}{2^{d-2}}\left(\tau^2-\left|\xi\right|^2-4s^2\right)^{\frac{d-2}{2}}\left(\tau^2-\left|\xi\right|^2\right)^{-\frac{1}{2}}
,
\end{equation*}
and hence
\begin{equation*}
I_{2,s}(\xi,\tau)=\frac{\left|\mathbb{S}^{d-1}\right|}{2^{\frac{d-1}{2}}}.
\end{equation*}
\end{lemma}
\proof As in \cite{Q'12} we use a one-parameter subgroup of
transformations $\left\{L^t\right\}_{t\in(-1,1)}$ of the group of
Lorentz transformations from $\sr^{d+1}$ to itself, defined as
\begin{equation*}
L^t(\xi,\tau)=\left(\frac{\xi_1 +t\tau}{\sqrt{1-t^2}}, \xi_2, \ldots, \xi_d,\frac{\tau+t\xi_1}{\sqrt{1-t^2}}\right),
\end{equation*}
for $(\xi, \tau)\in\sr^{d+1}$. We then note that the map
$(\xi,\tau)\rightarrow(A\xi,\tau)$ for any rotation $A$ of $\sr^d$
also belongs to the group of Lorentz transformations, and so for
fixed $(\xi,\tau)$ if we compose the operator $L^t$ where
$t=-\frac{\left|\xi\right|}{\tau}$ with the map described above
satisfying $A\xi=\left(\left|\xi\right|,0,\ldots,0\right)$ we obtain
a Lorentz transformation $\mathcal{L}$ so that
$\mathcal{L}(\xi,\tau)=(0,(\tau^2-\left|\xi\right|^2)^{\frac{1}{2}})$.
But then, as $\left|\det\mathcal{L}\right|=1$ it follows that the
convolution $\sigma_s\ast\sigma_s$ is also invariant under
$\mathcal{L}$, and hence
\begin{equation*}
J_s(\xi,\tau)=J_s(0, (\tau^2-\left|\xi\right|^2)^{\frac{1}{2}}), \;\; \tau>\left|\xi\right|.
\end{equation*}
This important reduction means that it suffices to consider
$J_s(0,z)$ for $z\in\sr$. Now
\begin{align*}
J_s(0,z) & =\int_{\sr^{d+1}}\delta(z-t-\phi_s(\left|-y\right|))\frac{1}{\phi_s(\left|-y\right|)}\delta(t-\phi_s(\left|y\right|))\frac{\,\mathrm{d}y}{\phi_s(\left|y\right|)}\,\mathrm{d}t \\
& =\int_{\sr^d}\delta(z-2\phi_s(\left|y\right|))\frac{\,\mathrm{d}y}{\phi_s(\left|y\right|)^2}.
\end{align*}
Using polar co-ordinates, we obtain
\begin{equation*}
\sigma_s\ast\sigma_s (0,z) = \left|\mathbb{S}^{d-1}\right|\int_0^\infty\delta(z-2\phi_s(r))\frac{r^{d-1}}{\phi_s(r)^2} \,\mathrm{d}r.
\end{equation*}
If we now make the change of variables $u=2\phi_s(r)$, by the
definition of $\phi_s$ we have that
$\frac{r}{\phi_s(r)^2}\mathrm{d}r=\frac{\mathrm{d}u}{u}$ and so
\begin{align*}
\sigma_s\ast\sigma_s (0,z) & =\left|\mathbb{S}^{d-1}\right|\int_{2s}^\infty\delta(z-u)\left(\frac{1}{2}\sqrt{u^2-4s^2}\right)^{d-2}\frac{\mathrm{d}u}{u} \\
& = \frac{\left|\mathbb{S}^{d-1}\right|}{2^{d-2}}\chi_{\left\{z\geq 2s\right\}}\frac{(z^2-4s^2)^{\frac{d-2}{2}}}{z},
\end{align*}
and the desired result follows from the Lorentz invariance discussed
above. \endproof

If we now integrate the inequality \eqref{KG-CS} for
$\left|\widetilde{u_2}\right|^2$ with respect to $\tau$ and $\xi$,
apply Plancherel's theorem and change the order of integration, we
obtain
\begin{align*}
\left\|u_2\right\|_{L_{x,t}^2}^2 &
=\frac{1}{(2\pi)^{d+1}}\left\|\smash{\widetilde{u_2}}\right\|_{L_{\xi,\tau}^2}^2 
\\
& \leq\frac{2^{-\frac{d-1}{2}}\left|\mathbb{S}^{d-1}\right|}{(2\pi)^{3d-1}}\int_{\sr^{2d}}\left|\smash{\widehat{f_1}}(y_1)\right|^2\left|\smash{\widehat{f_2}(y_2)}\right|^2 K_s(y)\phi_s(\left|y_1\right|)\phi_s(\left|y_2\right|)\,\mathrm{d}y.
\end{align*}
Moreover, if we consider the functions $f_j$ defined by
\begin{equation}\label{extremiser}
\phi_s(\left|y_j\right|)\widehat{f_j}(y_j)=e^{-a\phi_s(\left|y_j\right|)},
\end{equation}
for $a>0$ (and $j=1,2$), we immediately obtain that
\begin{equation*}
\widehat{F}(y)=\frac{e^{-a\tau}}{\sqrt{\phi_s(\left|y_1\right|)\phi_s(\left|y_2\right|)}},
\end{equation*}
on the support of the delta measures. Since the only place an
inequality was used was in the application of the Cauchy--Schwarz
inequality, this implies that we have equality for such functions.
Indeed, the above equality implies the existence of a scalar
function $g=g_s(\xi,\tau)$ so that
\begin{equation*}
K_s(y)\widehat{F}(y)=g(\xi,\tau)K_s(y)^{-1}(\phi_s(\left|y_1\right|)\phi_s(\left|y_2\right|))^{-\frac{1}{2}}
\end{equation*}
almost everywhere on the support of the delta measures, since on
this set $K_s$ may be written in terms of $\tau, \xi$ and $s$ only,
as shown above, and so may be absorbed into the function $g$. Hence
we have equality in \eqref{KG-CS} for these functions $f_j$, and
thus also in \eqref{bilinear-eq} for the constant
\begin{equation*}
\mathbf{KG}(d)=\frac{2^{-\frac{d-1}{2}}\left|\mathbb{S}^{d-1}\right|}{(2\pi)^{3d-1}},
\end{equation*}
implying that it is best possible.

\subsection{The case $d=1$ and $s>0$}
We note that formally, the calculation allowing us to derive
\eqref{bilinear-eq} also makes sense for $d=1$. However,
substituting $d=1$ into the expression for $K_s$ gives
\begin{align*}
K_s(y) & =\left[\left(\phi_s(\left|y_1\right|)\phi_s(\left|y_2\right|)-y_1y_2+s^2\right)\left(\phi_s(\left|y_1\right|)\phi_s(\left|y_2\right|)-y_1 y_2-s^2\right)\right]^{-\frac{1}{2}} \\
& =\left(\left(\phi_s(\left|y_1\right|)\phi_s(\left|y_2\right|)-y_1y_2\right)^2-s^4\right)^{-\frac{1}{2}} \\
& =\left(s^2(y_1^2+y_2^2)+2y_1^2y_2^2-2 y_1 y_2\phi_s(\left|y_1\right|)\phi_s(\left|y_2\right|)\right)^{-\frac{1}{2}},
\end{align*}
and since this weight is singular on the diagonal $\{(y_1,
y_2)\in\sr^2:y_1=y_2\}$, it is not difficult to construct a pair of
integrable functions $(f_1,f_2)$ for which the integral given by the
right hand side of \eqref{bilinear-eq} is unbounded. However if
$s>0$ the weight $K_s$ is well-defined for $y_1\neq y_2$ and if we
assume that $f_1$ and $f_2$ have disjointly supported Fourier
transforms, we have the identity \eqref{main-identity}. To prove
\eqref{main-identity} we follow a method used in \cite{F'70} for
restriction estimates on the sphere (see also \cite{H'73},
\cite{S'76}, \cite{C--S'72}, \cite{O--R'?}). Specifically, we write
\begin{align*}
e^{it\phi_s(\sqrt{-\Delta})}f_1(x) & \overline{e^{it\phi_s(\sqrt{-\Delta})}f_2(x)} \\
& = \int_{\sr^2}e^{ix(y_1-y_2)}e^{it((s^2+y_1^2)^{\frac{1}{2}}-(s^2+y_2^2)^{\frac{1}{2}})}\widehat{f_1}(y_1)\overline{\widehat{f_2}(y_2)}\,\mathrm{d}y_1\,\mathrm{d}y_2.
\end{align*}
If we make the change of variables $(y_1,y_2)\mapsto(u,v)$, where
$u=y_1-y_2$ and $v=\phi_s(\left| y_1 \right|)-\phi_s(\left| y_2
\right|)$, then the Jacobian will be
\begin{equation*}
\left|\det\left(
\begin{matrix}
1 & -1 \\
\frac{y_1}{\sqrt{s^2+y_1^2}} & -\frac{y_2}{\sqrt{s^2+y_2^2}}
\end{matrix}\right)
\right|^{-1}=\frac{\phi_s(\left| y_1 \right|)\phi_s(\left| y_2 \right|)}{\left|\phi_s(\left| y_1 \right|)y_2-\phi_s(\left| y_2 \right|)y_1\right|}.
\end{equation*}
Hence, we have
\begin{equation*}
e^{it\phi_s(\sqrt{-\Delta})}f_1(x) \overline{e^{it\phi_s(\sqrt{-\Delta})}f_2(x)}=\int_{\sr_+^2}e^{ixu}e^{itv}H(u,v)\,\mathrm{d}u\mathrm{d}v,
\end{equation*}
where $H$ is defined by
\begin{equation*}
H(u,v)=\frac{\phi_s(\left| y_1 \right|)\phi_s(\left| y_2 \right|)}{\left|\phi_s(\left| y_1 \right|)y_2-\phi_s(\left| y_2 \right|)y_1\right|}\widehat{f_1}(y_1)\widehat{f_2}(y_2).
\end{equation*}
By Plancherel's theorem,
\begin{equation*}
\left\|e^{it\phi_s(\sqrt{-\Delta})}f_1 e^{it\phi_s(\sqrt{-\Delta})}f_2\right\|_{L_{x,t}^2}^2=\frac{1}{(2\pi)^2}\left\|H\right\|_{L_{u,v}^2}^2=\frac{1}{(2\pi)^2}\int_{\sr^2}\left|H(u,v)\right|^2\,\mathrm{d}u\mathrm{d}v.
\end{equation*}
By reversing the change of variables done in the previous step, this
becomes
\begin{equation*}
\left\|u\right\|_{L^2}^2=\frac{1}{(2\pi)^2}\int_{\sr^2}\left|\smash{\widehat{f_1}(y_1)}\right|^2\left|\smash{\widehat{f_2}(y)}\right|^2\frac{\phi_s(\left| y_1 \right|)\phi_s(\left| y_2 \right|)}{\left|y_1\phi_s(\left| y_2 \right|)-y_2\phi_s(\left| y_1 \right|)\right|}\,\mathrm{d}y_1\,\mathrm{d}y_2.
\end{equation*}
Further, by a direct calculation, it is easily verified that 
\begin{align*}
\left((s^2+y_1^2)^{\frac{1}{2}}y_2-(s^2+y_2^2)^{\frac{1}{2}}y_1\right)^2 
= 
K_s(y)^{-2}.
\end{align*}
Note that from the above we can see that the only singularity of the
weight $K_s$ would be at a point in $\sr^2$ where
\begin{equation*}
y_1(s^2+y_2^2)^{\frac{1}{2}}=y_2(s^2+y_1^2)^{\frac{1}{2}},
\end{equation*}
which can only happen if $y_1=y_2$. It now remains to treat the case
$s=0$, where if we make no further assumptions than those used in
the case $s>0$, the argument breaks down. However if we assume that
$y_1y_2<0$ for all
$(y_1,y_2)\in\operatorname{supp}\widehat{f_1}\times\operatorname{supp}\widehat{f_2}$
(i.e.\ that the functions $f_1$ and $f_2$ on $\sr$ have disjoint
angular Fourier support) then it is not hard to see that the change
of variables analogous to $(y_1, y_2)\mapsto (u,v)$ makes sense and
the above argument yields an identity corresponding to
\eqref{main-identity} for $s=0$.

\section{Proof of Corollaries \ref{cor-2c} and \ref{cor-1c}} \label{section:corollaries}
We begin by establishing Corollary \ref{cor-1c} and show how to
deduce Corollary \ref{cor-2c}. Before proceeding, we recall the notation
\begin{equation*}
\left\|f\right\|_{(s)}^4=\left\|\phi_s(\sqrt{-\Delta})f\right\|_{L^2}^4-s^2\left\|(\phi_s(\sqrt{-\Delta}))^{\frac{1}{2}}f\right\|_{L^2}^4.
\end{equation*}
If we set $d=5$ and $f_1=f_2=f$ in \eqref{bilinear-eq}, then the
right hand side reduces to
\begin{align*}
\int_{\sr^{10}}& \left|\smash{\widehat{f}}(y_1)\right|^2 \left|\smash{\widehat{f}}(y_2)\right|^2  \phi_s(\left|y_1\right|)\phi_s(\left|y_2\right|)\frac{\left(\phi_s(\left|y_1\right|)\phi_s(\left|y_2\right|)-y_1\cdot y_2-s^2\right)^{\frac{3}{2}}}{\left(\phi_s(\left|y_1\right|)\phi_s(\left|y_2\right|)-y_1 \cdot y_2+s^2\right)^{\frac{1}{2}}}\,\mathrm{d}y \\
& \leq \int_{\sr^{10}}\left|\smash{\widehat{f}}(y_1)\right|^2 \left|\smash{\widehat{f}}(y_2)\right|^2 \phi_s(\left|y_1\right|)\phi_s(\left|y_2\right|)\left(\phi_s(\left|y_1\right|)\phi_s(\left|y_2\right|)-y_1 \cdot y_2-s^2\right)\,\mathrm{d}y \\
& = I_1-I_2,
\end{align*}
where
\begin{equation*}
I_1=\int_{\sr^{10}}\left(\left[\phi_s(\left|y_1\right|)\phi_s(\left|y_2\right|)\right]^2-s^2\phi_s(\left|y_1\right|)\phi_s(\left|y_2\right|)\right)\left|\smash{\widehat{f}}(y_1)\right|^2 \left|\smash{\widehat{f}}(y_2)\right|^2\,\mathrm{d}y_1 \,\mathrm{d}y_2,
\end{equation*}
and
\begin{equation*}
I_2=\int_{\sr^{10}}\left|\smash{\widehat{f}}(y_1)\right|^2 \left|\smash{\widehat{f}}(y_2)\right|^2 \phi_s(\left|y_1\right|)\phi_s(\left|y_2\right|)y_1\cdot y_2 \,\mathrm{d}y_1 \,\mathrm{d}y_2.
\end{equation*}
We can now use the observation of Carneiro in \cite{C'09},
\begin{equation}\label{C09}
\int_{\sr^{2d}}f(x)f(y)x\cdot y \,\mathrm{d}x\mathrm{d}y\geq 0
\end{equation}
which holds for any function $f$, with equality if $f$ is radial, to
obtain that $I_2\geq 0$. Hence, we have that
\begin{equation}\label{KG-1sided}
\left\|e^{it\phi_s(\sqrt{-\Delta})}f\right\|_{L_{x,t}^4}^4\leq (2\pi)^{10}\mathbf{KG}(5)\left\|f\right\|_{(s)}^4.
\end{equation}
Note however that we have used that
\begin{equation}\label{linfinity}
\frac{\phi_s(y_1)\phi_s(y_2)-y_1\cdot y_2-s^2}{\phi_s(y_1)\phi_s(y_2)-y_1\cdot y_2 +s^2}\leq 1,
\end{equation}
and this inequality is of course pointwise strict, but as with the $L^\infty$ analysis of the convolution of the measures $\sigma_s$ in \cite{Q'12} (Corollary 4.3, Lemma 4.4 and Lemma 4.5) we claim that when normalised, the functions $f_a$ form a maximising sequence for the inequality \eqref{KG-1sided}, 
as $a\rightarrow 0+$. As a consequence of this and inequality
\eqref{linfinity} we will obtain that the inequality
\eqref{KG-1sided} is sharp, and that there are no maximisers. We
recall that the functions $f_a$ are defined by
\begin{equation*}
\widehat{f_a}(x)=\frac{e^{-a\phi_s(\left|x\right|)}}{\phi_s(\left|x\right|)},
\end{equation*}
for $a>0$. By Theorem \ref{bilineart}, these satisfy inequality
\eqref{bilinear-eq} with equality, and by the observation after
inequality \eqref{C09} we also have that $I_2=0$ for such functions.
\begin{lemma}\label{e-s}
Suitably normalised, the functions $f_a$ form a maximising sequence for the
inequality \eqref{KG-1sided}. That is, we have that
\begin{equation*}
\operatorname*{lim}_{a\rightarrow 0+}\frac{\left\|\smash{e^{it\phi_s(\sqrt{-\Delta})}f_a}\right\|_{L_{x,t}^4}^4}{\left\|f_a\right\|_{(s)}^4}
=(2\pi)^{10}\mathbf{KG}(5).
\end{equation*}
\end{lemma}
\proof To prove Lemma \ref{e-s} we modify the approach in
\cite{Q'12}. Firstly, we calculate
\begin{align*}
(2\pi)^5\left\|\phi_s(\sqrt{-\Delta})^{\beta}f_a\right\|_{L^2}^2 & =\int_{\sr^5}\frac{e^{-2a\phi_s(\left|x\right|)}}{\phi_s(\left|x\right|)^{2-2\beta}}\,\mathrm{d}x \\
& =\left|\mathbb{S}^4\right|\int_0^\infty \frac{e^{-2a\sqrt{s^2+r^2}}}{(s^2+r^2)^{1-\beta}}r^4\,\mathrm{d}r \\
& =\left|\mathbb{S}^4\right|\int_s^\infty e^{-2au}\left(u^2-s^2\right)^{\frac{3}{2}}u^{2\beta-1}\,\mathrm{d}u \\
& =\frac{\left|\mathbb{S}^4\right|}{a^{2\beta}}\int_{as}^\infty e^{-2x}\left(\left(\frac{x}{a}\right)^2-s^2\right)^{\frac{3}{2}}x^{2\beta-1}\,\mathrm{d}x \\
& =\frac{\left|\mathbb{S}^4\right|}{a^{2\beta+3}}\int_{as}^\infty e^{-2x}\left(x^2-(as)^2\right)^{\frac{3}{2}}x^{2\beta-1}\,\mathrm{d}x
\end{align*}
for $\beta\in\left\{\frac{1}{2},1\right\}$, so that
\begin{equation}\label{h1-limit}
\operatorname*{lim}_{a\rightarrow 0+}a^5 (2\pi)^5 \left\|\phi_s(\sqrt{-\Delta})^\beta f_a\right\|_{L^2}^2=
\begin{cases}
\frac{3}{4}\left|\mathbb{S}^4\right| & \mbox{ if }\beta=1, \\
0 & \mbox{ if }\beta=\frac{1}{2}.
\end{cases}
\end{equation}
We now wish to evaluate
\begin{equation*}
\operatorname*{lim}_{a\rightarrow0+}a^{10}\left\|e^{it\phi_s(\sqrt{-\Delta})}f_a\right\|_{L^4}^4.
\end{equation*}
Observe that for these functions $f_a$ we can write this norm in
terms of the convolution of the measure $\sigma_s$ with itself.
Indeed, using Plancherel's theorem and then \eqref{conv-form} we
have
\begin{align*}
\left\|e^{it\phi_s(\sqrt{-\Delta})}f_a\right\|_{L^4}^4 & = \left\|\left(e^{it\phi_s(\sqrt{-\Delta})}f_a\right)^2\right\|_{L^2}^2 \\
& = \frac{1}{(2\pi)^{18}}\left\|\widetilde{e^{it\phi_s(\sqrt{-\Delta})}f_a}\ast\widetilde{e^{it\phi_s(\sqrt{-\Delta})}f_a}\right\|_{L^2}^2 \\
& = \frac{1}{(2\pi)^{14}}\left\|\int_{\sr^{10}} \frac{e^{-a(\phi_s(\left|y_1\right|\phi_s(\left|y_2\right|))}}{\phi_s(\left|y_1\right|)\phi_s(\left|y_2\right|)}\delta\binom{\tau-\phi_s(\left|y_1\right|)-\phi_s(\left|y_2\right|)}{\xi-y_1-y_2}\,\mathrm{d}x\right\|_{L_{\xi,\tau}^2}^2 \\
& = \frac{1}{(2\pi)^{14}}\left\|\int_{\sr^{10}} \frac{e^{-a\tau}}{\phi_s(\left|y_1\right|)\phi_s(\left|y_2\right|)}\delta\binom{\tau-\phi_s(\left|y_1\right|)-\phi_s(\left|y_2\right|)}{\xi-y_1-y_2}\,\mathrm{d}x\right\|_{L_{\xi,\tau}^2}^2 \\
& = \frac{1}{(2\pi)^{14}}\left\|e^{-a\tau}\sigma_s\ast\sigma_s(\tau, \xi)\right\|_{L_{\tau, \xi}^2}^2.
\end{align*}
By Lemma \ref{measure-conv}, we obtain
\begin{align*}
\left\|e^{it\phi_s(\sqrt{-\Delta})}f_a\right\|_{L^4}^4 & =\frac{\left|\mathbb{S}^4\right|^2}{2^6(2\pi)^{14}}\int_{\sr^{5+1}}e^{-2a\tau}\frac{(\tau^2-\left|\xi\right|^2-4s^2)^3}{\tau^2-\left|\xi\right|^2}\chi_{\left\{\tau\geq\sqrt{\left(2s\right)^2+\left|\xi\right|^2} \right\}}\,\mathrm{d}\xi \,\mathrm{d}\tau \\
& =\frac{\left|\mathbb{S}^4\right|^2}{2^6(2\pi)^{14}}\int_{\mathcal{H}_s} e^{-2a\tau}\left(\tau^2-\left|\xi\right|^2-4s^2\right)^2\left(1-\frac{4s^2}{\tau^2-\left|\xi\right|^2}\right)\,\mathrm{d}\xi \,\mathrm{d}\tau.
\end{align*}
To calculate the integral here, we write
\begin{align*}
\left(\tau^2-\left|\xi\right|^2-4s^2\right)^2 & =\left|\xi\right|^4+\tau^4-2\tau^2 \left|\xi\right|^2-8s^2\tau^2+8s^2\left|\xi\right|^2+16s^4 \\
& = \sum_{(j,k)\in\mathcal{T}}c_{j,k}\tau^{2j}\left|\xi\right|^{2k},
\end{align*}
where
\begin{equation*}
\displaystyle\mathcal{T}=\left\{\left(0,0\right), \left(1,0\right), \left(0,1\right), \left(1,1\right), \left(0,2\right), \left(2,0\right)\right\}\subseteq\mathbb{Z}\times\mathbb{Z}.
\end{equation*}
Hence, we obtain that
\begin{align*}
\int_{\mathcal{H}_s} e^{-2a\tau} & \left(\tau^2-\left|\xi\right|^2-4s^2\right)^2\left(1-\frac{4s^2}{\tau^2-\left|\xi\right|^2}\right)\,\mathrm{d}\xi \,\mathrm{d}\tau \\
& =\sum_{(j,k)\in\mathcal{T}}c_{j,k}\int_{2s}^\infty e^{-2a\tau}\tau^{2j}\int_{\left|\xi\right|\leq\sqrt{\tau^2-(2s)^2}}\left|\xi\right|^{2k}\left(1-\frac{4s^2}{\tau^2-\left|\xi\right|^2}\right)\, \mathrm{d}\xi\mathrm{d}\tau \\
& = \left|\mathbb{S}^4\right| \sum_{(j,k)\in\mathcal{T}}c_{j,k}\left(\mathrm{I}_{j,k}-4s^2 \mathrm{II}_{j,k}\right),
\end{align*}
where we define
\begin{equation*}
\mathrm{I}_{j,k}=\int_{2s}^\infty e^{-2a\tau}\tau^{2j}\int_0^{\sqrt{\tau^2-(2s)^2}}r^{2(k+2)}\,\mathrm{d}r\mathrm{d}\tau,
\end{equation*}
and
\begin{equation*}
\mathrm{II}_{j,k}=\int_{2s}^\infty e^{-2a\tau}\tau^{2j}\int_0^{\sqrt{\tau^2-(2s)^2}}\frac{r^{2(k+2)}}{\tau^2-r^2}\,\mathrm{d}r\mathrm{d}\tau.
\end{equation*}
\begin{claim}\label{claim1} We have that
\begin{equation*}
\operatorname*{lim}_{a\rightarrow 0+}a^{10}\mathrm{II}_{j,k}=\operatorname*{lim}_{a\rightarrow 0+}a^{10}\mathrm{II}_{j+1, k-1},
\end{equation*}
provided that $j+k<3$, $k>-\frac{3}{2}$ and $j\geq 0$, and
\begin{equation*}
\operatorname*{lim}_{a\rightarrow 0+}a^{10}\mathrm{II}_{j,-1}=0,
\end{equation*}
for $j\in\left\{1,2,3\right\}$.
\end{claim}
Assuming the claim to be true for the moment, it then follows that
$\mathrm{II}_{j,k}=0$ for each pair $(j,k)\in\mathcal{T}$, and hence
\begin{equation*}
\operatorname*{lim}_{a\rightarrow 0+}a^{10}\left\|e^{-a\tau}\sigma_s\ast\sigma_s(\tau, \xi)\right\|_{L_{\tau, \xi}^2}^2=\frac{\left|\mathbb{S}^4\right|^3}{2^6}\sum_{(j,k)\in\mathcal{T}}c_{j,k}\mathrm{I}_{j,k},
\end{equation*}
and we note that we can evaluate $\mathrm{I}_{j,k}$ directly, as
\begin{align*}
a^{10}\mathrm{I}_{j,k} & =\frac{a^{10}}{2k+5}\int_{2s}^\infty e^{-2a\tau}\tau^{2j}(\tau^2-(2s)^2)^{\frac{2k+5}{2}}\,\mathrm{d}\tau \\
& =\frac{1}{2k+5}a^{4-2(j+k)}\int_{2as}^\infty e^{-2x}x^{2j}(x^2-(2as)^2)^{\frac{2k+5}{2}}\,\mathrm{d}x.
\end{align*}
Hence, since the latter integral converges for any $a>0$ we have
\begin{equation*}
\operatorname*{lim}_{a\rightarrow 0+}a^{10}\mathrm{I}_{j,k}=
\begin{cases}
\displaystyle\frac{1}{2k+5}\int_0^\infty e^{-2x}x^9\,\mathrm{d}x & \mbox{ if } j+k=2, \\
0 & \mbox{ if } j+k<2.
\end{cases}
\end{equation*}
In all, since
\begin{equation*}
\int_0^\infty x^{\ell}e^{-2x}\,\mathrm{d}x=\frac{\ell !}{2^{\ell+1}}
\end{equation*}
we obtain that
\begin{equation*}
\operatorname*{lim}_{a\rightarrow 0+}a^{10}\left\|e^{it\phi_s(\sqrt{-\Delta})}f_a\right\|_{L^4}^4=\frac{\left|\mathbb{S}^4\right|^{3}}{(2\pi)^{14}}\left(\frac{1}{5}+\frac{1}{9}-\frac{2}{7}\right)\frac{9!}{2^{16}},
\end{equation*}
and
\begin{equation*}
\operatorname*{lim}_{a\rightarrow 0+}a^{10}(2\pi)^{10}\left\|f_a\right\|_{(s)}^4
=\frac{9}{16}\left|\mathbb{S}^4\right|^2.
\end{equation*}
Hence
\begin{equation*}
\operatorname*{lim}_{a\rightarrow 0+}\frac{\left\|e^{it\phi_s(\sqrt{-\Delta})}f_a\right\|_{L^4}^4}{(2\pi)^{10}\left\|f_a\right\|_{(s)}^4}
= \mathbf{KG}(5),
\end{equation*}
as claimed, and therefore the constant $(2\pi)^{10}\mathbf{KG}(5)$
is best possible for the inequality $\eqref{cor-1.1}$. We also
remark at this point that the constant $(2\pi)^{10}\mathbf{KG}(5)$
is also sharp for the inequality \eqref{l4h1}; as discussed in
Section 1 this follows from the sequence $\left(g_a\right)_{a>0}$
defined by \eqref{g_a}, since by \eqref{h1-limit} we have
$a^5\left\|\phi_s(\sqrt{-\Delta})^{\frac{1}{2}}f_a\right\|^2_{L^2}\rightarrow
0$ as $a\rightarrow 0+$. It now remains to prove Claim \ref{claim1}.


\emph{Proof of Claim \ref{claim1}}. For the first part, we note that
\begin{align*}
\int_0^{\sqrt{\tau^2-(2s)^2}}\frac{r^{2(k+2)}}{\tau^2-r^2}\,\mathrm{d}r & =\int_0^{\sqrt{\tau^2-(2s)^2}}r^{2(k+1)}\left(\frac{\tau^2}{\tau^2-r^2}-1\right)\,\mathrm{d}r \\
& =\tau^2 \int_0^{\sqrt{\tau^2-(2s)^2}}\frac{r^{2(k+1)}}{\tau^2-r^2}\,\mathrm{d}r-\frac{1}{2k+3}\left(\tau^2-(2s)^2\right)^{k+\frac{3}{2}}.
\end{align*}
But then,
\begin{align*}
a^{10}\int_{2s}^\infty e^{-2a\tau}\tau^{2j}\left(\tau^2-(2s)^2\right)^{k+\frac{3}{2}}\,\mathrm{d}\tau & =a^{6-2(j+k)}\int_{2as}^\infty e^{-2x}x^{2j}\left(x^2-(2as)^2\right)^{k+\frac{3}{2}}\\
& \rightarrow 0
\end{align*}
as $a\rightarrow 0+$ since $j+k<3$ and $j\geq 0$. Hence,
\begin{align*}
\operatorname*{lim}_{a\rightarrow 0+}a^{10}\mathrm{II}_{j,k}& =\operatorname*{lim}_{a\rightarrow 0+}a^{10}\int_{2s}^\infty e^{-2a\tau}\tau^{2(j+1)}\int_0^{\sqrt{\tau^2-(2s)^2}}\frac{r^{2(k+1)}}{\tau^2-r^2}\,\mathrm{d}r\mathrm{d}\tau \\
& =\operatorname*{lim}_{a\rightarrow 0+}a^{10}\mathrm{II}_{j+1, k-1}.
\end{align*}
For the second part, by a simple change of variables we can
calculate $\mathrm{I}_{j,-1}$ directly, as in \cite{Q'12}. We have
\begin{align*}
a^{10}\mathrm{II}_{j,-1} & =a^{8-2j}\int_{2as}^\infty e^{-2x}x^{2j+1}\operatorname{log}\left(\frac{x+\sqrt{x^2+(2as)}}{2as}\right)\,\mathrm{d}x \\
& = a^{8-2j}\int_{2as}^\infty e^{-2x}x^{2j+1}\operatorname{log}\left(x+\sqrt{x^2+(2as)}\right)\,\mathrm{d}x \\
& \;\;\;\; -a^{8-2j}\operatorname{log}(2as)\int_{2as}^\infty e^{-2x}x^{2j+1}\,\mathrm{d}x \\
& \rightarrow 0
\end{align*}
as $a\rightarrow 0+$, since $8-2j>0$ and $2j+1>0$.
\endproof

We conclude the section by showing how Corollary \ref{cor-2c} is
deduced from Corollary \ref{cor-1c}; to do this we follow the
approach of Foschi in \cite{F'06}. Suppose $u$ solves
\begin{equation}\label{KG-eq}
\partial_{tt} u-\Delta u+ u=0,
\end{equation}
where $u=u(x,t)$ is a function defined on $\sr^{5+1}$. Then we can
write $u=u_+ +u_-$, where
\begin{equation}\label{uplusuminus}
u_+=e^{it\phi_1(\sqrt{-\Delta})}f_+, \;\;\;\; u_-=e^{-it\phi_1(\sqrt{-\Delta})}f_-,
\end{equation}
for functions $f_+$ and $f_-$ defined using the initial data by
\begin{equation*}
u(0)=f_+ +f_-, \;\;\;\; \partial_t u(0)=i\phi_1(\sqrt{-\Delta})(f_+ -f_-).
\end{equation*}
Then
\begin{equation*}
\left\|u\right\|_{L^4}^4 =\left\|u_+ +u_-\right\|_{L^4}^4=\left\|u_+^2+u_-^2+2u_+u_-\right\|_{L^2}^2.
\end{equation*}
We claim that the supports of the space-time Fourier transforms of
the functions $u_+^2$, $u_-^2$ and $u_+ u_-$ are pairwise disjoint.
We have already seen that
\begin{equation*}
\operatorname{supp}\widetilde{u_+}\subseteq\left\{\left(\xi,\tau\right)\in\sr^{d+1}:\tau\geq\sqrt{4+\left|\xi\right|^2}\right\},
\end{equation*}
and by an identical argument we have that
\begin{equation*}
\operatorname{supp}\widetilde{u_-}\subseteq\left\{\left(\xi,\tau\right)\in\sr^{d+1}:\tau\leq-\sqrt{4+\left|\xi\right|^2}\right\}.
\end{equation*}
It remains to show that
\begin{equation*}
\operatorname{supp}\widetilde{u_+u_-}\subseteq\left\{\left(\xi,\tau\right)\in\sr^{d+1}:\left|\tau\right|\leq\sqrt{4+\left|\xi\right|^2}\right\}.
\end{equation*}
We note that this was shown in \cite{Q'12}, we include it here for
completeness. To see this, note that analogously to
\eqref{conv-form} we will have, for $(\xi,\tau)\in\sr^{d+1}$,
\begin{equation*}
\widetilde{u_+u_-}(\xi,\tau)=\frac{1}{(2\pi)^{d-1}}\int_{\sr^{2d}}\frac{\widehat{F}(y)}{(1+\left|y_1\right|^2)^{\frac{1}{4}}(1+\left|y_2\right|^2)^{\frac{1}{4}}}\delta\binom{\tau-\phi_1(\left|y_1\right|)+\phi_1(\left|y_2\right|)}{\xi-y_1-y_2}\,\mathrm{d}y.
\end{equation*}
Set $\xi=y_1+y_2$ and $\tau=\phi_1(\left| y_1 \right|)-\phi_1(\left|
y_2 \right|)$, then we have
\begin{equation*}
\left|\xi\right|^2=\left|y_1\right|^2+\left|y_2\right|^2+2y_1\cdot y_2,
\end{equation*}
and,
\begin{equation*}
\tau^2=2s^2+\left|y_1\right|^2+\left|y_2\right|^2-2\phi_1(\left| y_1 \right|)\phi_1(\left| y_2 \right|),
\end{equation*}
so that
\begin{equation*}
\tau^2-\left|\xi\right|^2=2-2\phi_1(\left| y_1 \right|)\phi_1(\left| y_2 \right|)-2y_1\cdot y_2\leq 4,
\end{equation*}
where the final inequality follows from \eqref{ineq-1}. Hence,
\begin{equation*}
\left\|u\right\|_{L^4}^4=\left\|u_+^2\right\|_2^2+\left\|u_-^2\right\|_2^2+4\left\|u_+u_-\right\|_2^2.
\end{equation*}
Combining the above equality with the sharp polynomial inequality
for non-negative real numbers $X$ and $Y$
\begin{equation*}
X^2+Y^2+4XY\leq \frac{3}{2}(X+Y)^2
\end{equation*}
with equality if and only if $X=Y$, we obtain the following sharp
inequality in terms of the one-sided propagators:
\begin{equation*}
\left\|u\right\|_{L^4}^4 \leq \frac{3}{2}\left(\left\|u_+\right\|_{L^4}^2+\left\|u_-\right\|_{L^4}^2\right)^2.
\end{equation*}
Using the inequality \eqref{l4h1}, it follows that
\[\left\|u\right\|_{L^4}^4\leq \frac{1}{16\pi^2}\left(\left\|f_+\right\|_{H^1}^2+\left\|f_-\right\|_{H^1}^2\right)^2.\]
But then, by the definition of $f_+$ and $f_-$ and the parallelogram
law, the right hand side equals
\[\frac{1}{16\pi^2}\left(\frac{1}{2}\left\|u(0)\right\|_{H^1}^2+\frac{1}{2}\left\|\phi_1(\sqrt{-\Delta})^{-1}\partial_tu(0)\right\|_{H^1}^2\right)^2=\frac{1}{64\pi^2}\left(\left\|u(0)\right\|_{H^1}^2+\left\|\partial_tu(0)\right\|_{L^2}^2\right)^2,\]
which completes the proof of Corollary \ref{cor-2c}.

\qed
\section{Proof of Proposition 1}
In this section, it will be convenient to abuse notation slightly
and think of $\phi_s(x)$ as a function on $\sr^5$ by identifying
with $\phi_s(\left|x\right|)$. If we let
\begin{equation*}
\left\|g \otimes g\right\|_{(\tau,\xi)}^2=\int_{\sr^{10}}\left|g(y_1)\right|^2 \left|g(y_2)\right|^2\delta\binom{\tau-\phi_s\left(\left|y_1\right|\right)-\phi_s(\left|y_2\right|)}{\xi-y_1-y_2}\,\mathrm{d}y_1\mathrm{d}y_2,
\end{equation*}
by the proof of the bilinear inequality of Theorem \ref{bilineart}
and by Lemma \ref{measure-conv}, we have that
\begin{align*}
& \left\|e^{it\phi_s(\sqrt{-\Delta})}g_n\right\|_{L^4(\sr^6)}^4 \\
& \qquad \qquad \leq \frac{1}{24\pi^2}\int_{\mathcal{H}_s}\left\|\smash{\phi_s^{\frac{1}{2}}\widehat{g_n}} \otimes  \smash{\phi_s^{\frac{1}{2}}\widehat{g_n}}  \right\|_{(\tau,\xi)}^2 \frac{(\tau^2-\left|\xi^2\right|-4s^2)^{\frac{3}{2}}}{(\tau^2-\left|\xi\right|^2)^{\frac{1}{2}}}\,\mathrm{d}\xi\mathrm{d}\tau \\
& \qquad \qquad = \frac{1}{24\pi^2}\int_{\mathcal{H}_s}\left\|\smash{\phi_s^{\frac{1}{2}}\widehat{g_n}} \otimes \smash{\phi_s^{\frac{1}{2}}\widehat{g_n}} \right\|_{(\tau,\xi)}^2\left(\tau^2-\left|\xi\right|^2\right)\left(1-\frac{4s^2}{\tau^2-\left|\xi\right|^2}\right)^{\frac{3}{2}}\,\mathrm{d}\tau\mathrm{d}\xi \\
& \qquad \qquad = \frac{1}{24\pi^2}\left(\left\|g_n\right\|_{(s)}^4-\mathcal{I}_n-\mathcal{J}_n\right) \\
& \qquad \qquad \leq \frac{1}{24\pi^2}\left(1-\mathcal{I}_n-\mathcal{J}_n\right),
\end{align*}
where
\begin{equation*}
\mathcal{I}_n=\int_{\sr^{10}}\left|\smash{\widehat{g_n}(y_1)}\right|^2 \left|\smash{\widehat{g_n}(y_2)}\right|^2 \phi_s(\left|y_1\right|)\phi_s(\left|y_2\right|)y_1 \cdot y_2\,\mathrm{d}y_1 \mathrm{d}y_2\geq 0,
\end{equation*}
and
\begin{equation*}
\mathcal{J}_n=\int_{\mathcal{H}_s}\left\|\smash{\phi_s^{\frac{1}{2}}\widehat{g_n}} \otimes \smash{\phi_s^{\frac{1}{2}}\widehat{g_n}}\right\|_{(\tau,\xi)}^2\left(\tau^2-\left|\xi\right|^2\right)\left(1-\left(1-\frac{4s^2}{\tau^2-\left|\xi\right|^2}\right)^{\frac{3}{2}}\right)\,\mathrm{d}\tau\mathrm{d}\xi\geq 0.
\end{equation*}

But then since $(g_n)_{n\geq 1}$ is a maximising sequence for
inequality \eqref{cor-1.1}, it follows that $\mathcal{I}_n,
\mathcal{J}_n\rightarrow 0$ as
$n\rightarrow\infty$. 
Firstly, since $0<1-\frac{4s^2}{\tau^2-\left|\xi\right|^2}<1$, 
we obtain
\begin{equation*}
\left(\int_{\sr^5}\left|\smash{\widehat{g_n}(y_1)}\right|^2\phi_s(\left|y_1\right|)\,\mathrm{d}y_1\right)^2< C\mathcal{J}_n,
\end{equation*}
for some positive constant $C$, and hence
\begin{equation*}
\int_{\sr^5}\left|\smash{\widehat{g_n}(y_1)}\right|^2\phi_s(\left|y_1\right|)\,\mathrm{d}y_1 \rightarrow 0
\end{equation*}
as $n\rightarrow\infty$.

Now, to prove \eqref{conc-ineq}, using the fact that on the delta measures we have that
\begin{equation}\label{id-2}
\phi_s(\left|y_1\right|)\phi_s(\left|y_2\right|)-y_1\cdot y_2+s^2=\frac{1}{2}\left(\tau^2-\left|\xi\right|^2\right),
\end{equation}
if $y_1, y_2\in B(0,R)$ it is easy to see that for such $\tau, \xi$,
\begin{equation*}
\tau^2-\left|\xi\right|^2\leq 2(R^2+s^2).
\end{equation*}
Thus,
\begin{align*}
& \int_{\mathcal{H}^s}\int_{B(0,R)}\left|\smash{\widehat{g_n}(y_1)} \right|^2 \left|\smash{\widehat{g_n}(y_2)}\right|^2\phi_s(\left|y_1\right|)\phi_s(\left|y_2\right|)\left(\tau^2-\left|\xi\right|^2\right)\\
&  \qquad \qquad \qquad \qquad \qquad \qquad \times
\delta\binom{\tau-\phi_s(\left|y_1\right|)-\phi_s(\left|y_2\right|)}{\xi-y_1-y_2}\mathrm{d}y\mathrm{d}\xi\mathrm{d}\tau \\
& \qquad \qquad \qquad \leq 2(R^2+s^2)\left(\int_{\sr^5}\left|\smash{\widehat{g_n}(y_1)}\right|^2\phi_s(\left|y_1\right|)\mathrm{d}y_1\right)^2\rightarrow 0
\end{align*}
as $n\rightarrow\infty$. Using \eqref{id-2} we obtain
\begin{align*}
& \left(\int_{B(0,R)}\left|\smash{\widehat{g_n}}(y_1)\right|^2\phi_s(\left|y_1\right|)^2\mathrm{d}y_1\right)^2 \\
& \leq
\int_{B(0,R)}\left|\widehat{g_n}(y_1)\right|^2\left|\smash{\widehat{g_n}}(y_2)\right|^2\left(\phi_s(\left|y_1\right|)^2\phi_s(\left|y_2\right|)^2+s^2\phi_s(\left|y_1\right|)\phi_s(\left|y_2\right|)\right)\mathrm{d}y_1\mathrm{d}y_2 \\
& =\int_{B(0,R)}\left|\smash{\widehat{g_n}}(y_1)\right|^2\left|\smash{\widehat{g_n}}(y_2)\right|^2 \phi_s(\left|y_1\right|)\phi_s(\left|y_2\right|)y_1 \cdot y_2 \,\mathrm{d}y_1\mathrm{d}y_2 \\
& \qquad + \frac{1}{2}\int_{\mathcal{H}^s}\int_{B(0,R)}\left|\smash{\widehat{g_n}(y_1)} \right|^2 \left|\smash{\widehat{g_n}(y_2)}\right|^2\phi_s(\left|y_1\right|)\phi_s(\left|y_2\right|)\left(\tau^2-\left|\xi\right|^2\right)\\
& \qquad \qquad \qquad \qquad \qquad \qquad \times
\delta\binom{\tau-\phi_s(\left|y_1\right|)-\phi_s(\left|y_2\right|)}{\xi-y_1-y_2}\,\mathrm{d}y\mathrm{d}\xi\mathrm{d}\tau \\
& \leq \mathcal{I}_n+(R^2+s^2)\left(\int_{\sr^5}\left|\smash{\widehat{g_n}(y_1)}\right|^2\phi_s(\left|y_1\right|)\,\mathrm{d}y_1\right)^2\rightarrow 0
\end{align*}
as $n\rightarrow\infty$. But then if $\varepsilon, R$ are given then
we can choose $N$, as desired.
\endproof
\begin{acknowledgement}
I wish to express my gratitude to my thesis advisor, Dr. N. Bez, for introducing me to this problem and for many interesting discussions.
\end{acknowledgement}

\end{document}